\documentclass{amsproc}
\usepackage{amssymb}
\usepackage{amsmath}
\usepackage[matrix,arrow,curve,cmtip]{xy}

\numberwithin{equation}{section}

\theoremstyle{plain}   
\newtheorem{bigthm}{Theorem}   

\newtheorem{theorem}[equation]{Theorem}  
     
\newtheorem{lemma}[equation]{Lemma}         
\newtheorem{prop}[equation]{Proposition}

\theoremstyle{definition}

\theoremstyle{remark}
\newtheorem{remark}[equation]{Remark}

\newcommand{\TF}{\operatorname{TF}}
\newcommand{\TC}{\operatorname{TC}}
\newcommand{\TR}{\operatorname{TR}}

\newcommand{\HH}{\operatorname{HH}}
\newcommand{\HC}{\operatorname{HC}}

\newcommand{\N}{\mathbb{N}}
\newcommand{\Z}{\mathbb{Z}}
\newcommand{\Q}{\mathbb{Q}}
\newcommand{\R}{\mathbb{R}}
\newcommand{\C}{\mathbb{C}}

\newcommand{\Fp}{\mathbb{F}_p}

\newcommand{\cy}{{\text{\rm cy}}}

\newcommand{\id}{\operatorname{id}}

\newcommand{\xto}{\xrightarrow}

\begin{document}

\title{On the $K$-theory of the coordinate axes in the plane}

\author{Lars Hesselholt}

\address{Massachusetts Institute of Technology, Cambridge,
Massachusetts}

\email{larsh@math.mit.edu}

\address{Nagoya University, Nagoya, Japan}

\email{larsh@math.nagoya-u.ac.jp}

\thanks{The author was partially supported by COE (Japan) and the
  National Science Foundation}

\begin{abstract}
Let $k$ be a regular $\Fp$-algebra, let $A = k[x,y]/(xy)$ be the
coordinate ring of the coordinate axes in the affine $k$-plane, and
let $I = (x,y)$ be the ideal that defines the intersection point. We
evaluate the relative $K$-groups $K_q(A,I)$ in terms of the
groups of big de~Rham-Witt forms of $k$. The result generalizes
previous results for $K_1$ and $K_2$ by Dennis and Krusemeyer.
\end{abstract}

\maketitle

\section*{Introduction}

Let $k$ be a ring, and let $A = k[x,y]/(xy)$ be the coordinate ring of
the coordinate axes in the affine $k$-plane. The $K$-groups of $A$
decompose as the direct sum
$$K_q(A) = K_q(k) \oplus K_q(A,I)$$
of the $K$-groups of the ground ring $k$ and the relative $K$-groups
of $A$ with respect to the ideal $I = (x,y)$. In this paper we
evaluate the groups $K_q(A,I)$ completely in the case where $k$ is a
regular $\Fp$-algebra. The result is stated in terms of the groups of
big de~Rham-Witt forms of $k$ as follows.

\begin{bigthm}\label{main}Let $k$ be a regular\, $\Fp$-algebra, let $A
= k[x,y]/(xy)$ be the coordinate ring of the coordinate axes in the
affine $k$-plane, and let $I \subset A$ be the ideal generated by $x$
and $y$. Then for all integers $q$, there is a canonical isomorphism
$$K_q(A,I) \xleftarrow{\sim} \bigoplus_{m \geqslant 1}
\mathbf{W}_m\Omega_k^{q - 2m}$$
where $\mathbf{W}_m\Omega_k^j$ is the group of big de~Rham-Witt
$j$-forms of $k$.
\end{bigthm}

The group $K_2(A,I)$ was evaluated by Dennis and
Krusemeyer~\cite{denniskrusemeyer} twenty-five years ago. But it was
previously known only that the higher relative $K$-groups are
$p$-primary torsion groups~\cite{weibel2}. The group of big
de~Rham-Witt $j$-forms $\mathbf{W}_n\Omega_k^j$ was introduced
in~\cite[Def.~1.1.6]{hm2}. It decomposes as a product of the more familiar
$p$-typical de~Rham-Witt $j$-forms $W_s\Omega_k^j$ defined by
Bloch-Deligne-Illusie~\cite{illusie}. Indeed,
by~\cite[~Cor.~1.2.6]{hm2} there is a canonical isomorphism
$$\mathbf{W}_m\Omega_k^j \xto{\sim} \prod_{d} W_s\Omega_k^j,$$
where the product ranges over all integers $1 \leqslant d \leqslant m$ that are
not divisible by $p$, and where $s = s(m,d)$ is the unique positive
integer with $p^{s-1}d \leqslant m < p^sd$. The structure of the groups
$W_s\Omega_k^j$ is well-understood by~\cite[I.3.9]{illusie}. For
example, $W_s\Omega_{\Fp}^j$ is canonically isomorphic
to $\Z/p^s\Z$, for $j = 0$, and is equal to zero, for $j > 0$.

Let $B = k[x] \times k[y]$ be the normalization of the ring $A$, and
let $K(A,B,I)$ be the bi-relative $K$-theory spectrum defined to be the
iterated mapping fiber of the following diagram of $K$-theory
spectra.
$$\xymatrix{
{ K(A) } \ar[r] \ar[d] &
{ K(A/I) } \ar[d] \cr
{ K(B) } \ar[r] &
{ K(B/I) } \cr
}$$
The lower horizontal map in this diagram is a weak equivalence since
the ring $k$ is regular. It follows that also the canonical map
$$K(A,B,I) \to K(A,I)$$
is a weak equivalence. We proved recently in~\cite[Thm.~A]{gh4} that
for every prime $p$, the cyclotomic trace map induces an
isomorphism
$$K_q(A,B,I,\Z/p^v) \xto{\sim} \TC_q(A,B,I;p,\Z/p^v)$$
and it is the bi-relative topological cyclic homology groups on the
right-hand side that we evaluate here. The method is similar to the
calculation of the topological cyclic homology of the ring of dual
numbers by the author and Madsen~\cite{hm,hm2}. We first prove a
general formula that expresses the
bi-relative topological cyclic homology groups above in terms of the
$RO(\mathbb{T})$-graded equivariant homotopy groups
$$\TR_{\alpha}^n(k;p) = [S^{\alpha} \wedge
(\mathbb{T}/C_{p^{n-1}})_+, T(k)]_{\mathbb{T}}$$
of the topological Hochschild $\mathbb{T}$-spectrum $T(k)$. Here
$\mathbb{T}$ is the circle group and $C_r \subset \mathbb{T}$ is the
subgroup of order $r$. To state the formula, which is valid for any
ring $k$, we let $\lambda_i$ be the complex
$\mathbb{T}$-representation $\C(1) \oplus \dots \oplus \C(i)$.

\begin{bigthm}\label{tcthm}Let $k$ be an\, $\Fp$-algebra, let $A =
k[x,y]/(xy)$ and $B = k[x] \times k[y]$, and let $I$ be the common
ideal of $A$ and $B$ generated by $x$ and $y$. Then for all integers
$q$, there is a canonical isomorphism
$$\TC_q(A,B,I;p) \xrightarrow{\sim} \prod \lim_R
\TR_{q-\lambda_{p^{r-1}d}}^r(k;p)$$
where the product ranges over the positive integers $d$ that are not
divisible by $p$. The analogous statement for the groups with
$\Z/p^v$-coefficients is valid for any ring $k$. 
\end{bigthm}

The limit system on the right-hand side of the statement of
Thm.~\ref{tcthm} stabilizes in the sense that the structure map
$$R \colon \TR_{q-\lambda_{p^{r-1}d}}^r(k;p) \to
\TR_{q-\lambda_{p^{r-2}d}}^{r-1}(k;p)$$
is an isomorphism for $q < \dim_{\R}(\lambda_{p^{r-1}d})$. See
Lemma~\ref{mittagleffler} below.

If $k$ is a regular $\Fp$-algebra, the structure of the groups on the
right-hand side of the statement of Thm.~\ref{tcthm} was determined 
in~\cite[Thm.~2.2.2]{hm2}; see also~\cite[Thm.~16]{h3}. We recall the
result in Sect.~3 below and complete the proof of Thm.~\ref{main}.

Finally, we mention that the analog of Thm.~\ref{main} for $k$ a
regular $\Q$-algebra is known. Indeed, by a recent theorem of
Corti\~{n}as~\cite[Thm.~0.1]{cortinas} (which inspired us to
prove~\cite[Thm.~A]{gh4}), the trace map induces an isomorphism
$$K_q(A,B,I) \otimes \Q \xto{\sim} \HC_q^{-}(A \otimes \Q, B \otimes
\Q , I \otimes \Q)$$
and the bi-relative negative cyclic homology groups on the right-hand
side were evaluated long ago by Geller, Reid, and
Weibel~\cite{gellerreidweibel,gellerreidweibel1}. The result is that,
if $k$ is a regular $\Q$-algebra, then there is a canonical
isomorphism
$$K_q(A,I) \xleftarrow{\sim} \bigoplus_{m \geqslant 1} \Omega_k^{q-2m}$$
where $\Omega_k^j$ is the group of absolute K\"ahler $j$-forms of
$k$. This formula differs from the formula of Thm.~\ref{main} in
degrees $q \geqslant 4$. Indeed, the group $K_4(A,I)$ is isomorphic to
$\Omega_k^2 \oplus k$, if $k$ is a regular $\mathbb{Q}$-algebra, to
$\Omega_k^2 \oplus W_2(k)$, if $k$ is a regular
$\mathbb{F}_2$-algebra, and to $\Omega_k^2 \oplus k \oplus k$, if $k$
is a regular $\mathbb{F}_p$-algebra and $p > 2$. We also remark that,
if $k$ is a regular $\mathbb{F}_p$-algebra, and if $p^{s-1} \leqslant
n < p^s$, then the group $K_{2n}(A,I)$ has exponent exactly $p^s$.

The result of Thm.~\ref{main} was announced in~\cite[Thm.~C]{gh4}.

All rings considered in this paper are assumed to be commutative. We
write $\mathbb{N}$ and $I_p$ for the sets of positive integers and
positive integers prime to $p$, respectively. We say that a map of
$\mathbb{T}$-spectra is an $\mathcal{F}$-equivalence if the induced
map of $C$-fixed point spectra is a weak equivalence, for all finite
subgroups $C \subset \mathbb{T}$.

Finally, the author would like to thank an anonymous referee for a
very careful reading of an earlier version of this paper and for a
number of helpful suggestions on improving the exposition.

\section{Topological Hochschild homology}

The proof of Thm.~\ref{tcthm} of the introduction is based on a
description of the bi-relative topological Hochschild $\mathbb{T}$-spectrum
$T(A,B,I)$ defined to be the iterated mapping fiber of the following
diagram of topological Hochschild $\mathbb{T}$-spectra.
\begin{equation}\label{birelative}
\xymatrix{
{ T(A) } \ar[r] \ar[d] &
{ T(A/I) } \ar[d] \cr
{ T(B) } \ar[r] &
{ T(B/I). } \cr
}
\end{equation}
We refer~\cite{h3} for an introduction to topological
Hochschild and cyclic homology and for further references. In this
section, we shall prove the following result. 

\begin{prop}\label{thh}Let $k$ be any ring, let $A = k[x,y]/(xy)$ and
$B = k[x] \times k[y]$, and let $I$ be the common ideal generated by
$x$ and $y$. Then there is a canonical $\mathcal{F}$-equivalence of
  $\mathbb{T}$-spectra
$$\bigvee_{i \in \N} T(k) \wedge S^{\lambda_i} \wedge
(\mathbb{T}/C_i)_+[1] \xto{\sim} T(A,B,I)$$
where, on the left-hand side, $[1]$ indicates desuspension.
\end{prop}

The rings $B = k[x] \times k[y]$ and $B/I = k \times k$ are both
product rings. Moreover, topological Hochschild homology preserves
products of rings in the sense that for every pair of rings $R$ and
$S$, the canonical map of $\mathbb{T}$-spectra
$$T(R \times S) \to T(R) \times T(S)$$
is an $\mathcal{F}$-equivalence~\cite[Prop.~4.20]{bokstedthsiangmadsen}. 
Hence, the canonical map from $T(A,B,I)$ to the iterated mapping fiber
of the following diagram of $\mathbb{T}$-spectra is an
$\mathcal{F}$-equivalence.
\begin{equation}\label{square}
\xymatrix{
{ T(k[x,y]/(xy)) } \ar[r]^(.61){\epsilon} \ar[d]^{(\phi,\phi')} &
{ T(k) } \ar[d]^{\Delta} \cr
{ T(k[x]) \times T(k[y]) } \ar[r]^(.55){\epsilon\times\epsilon} &
{ T(k) \times T(k). } \cr
}
\end{equation}
The rings that occur in this diagram are all pointed monoid algebras. 
By definition, a pointed monoid $\Pi$ is a monoid in the symmetric
monoidal category of pointed sets and smash product, and the pointed
monoid algebra $k(\Pi)$ is the quotient of the monoid algebra $k[\Pi]$
by the ideal generated by the base-point of $\Pi$. The diagram of
rings~(\ref{square}) is then induced from the diagram of pointed
monoids
$$\xymatrix{
{ \Pi^2 } \ar[r]^{\epsilon} \ar[d]^{(\phi,\phi')} &
{ \Pi^0 } \ar[d]^{\Delta} \cr
{ \Pi^1 \times \Pi^1 } \ar[r]^{\epsilon \times \epsilon} &
{ \Pi^0 \times \Pi^0 } \cr
}$$
where $\Pi^0 = \{0, 1\}$ with base-point $0$, where $\Pi^1 = \{0, 1,
z, z^2, \dots \}$ with base-point $0$, and where $\Pi^2 = \{0, 1, x,
x^2, \dots, y, y^2, \dots \}$ with base-point $0$ and with
multiplication given by $xy = 0$. The map $\phi$ (resp.~$\phi'$) takes
the variables $x$ and $y$ to $z$ and $0$ (resp.~ to $0$ and $z$), and
the maps labeled $\epsilon$ take the variables $x$, $y$, and $z$ to
$1$.

The topological Hochschild $\mathbb{T}$-spectrum of a pointed monoid
algebra $k(\Pi)$ decomposes, up to $\mathcal{F}$-equivalence, as the
smash product
\begin{equation}\label{separate}
T(k) \wedge N^{\cy}(\Pi) \xto{\sim} T(k(\Pi))
\end{equation}
of the topological Hochschild $\mathbb{T}$-spectrum of the coefficient
ring $k$ and the cyclic bar-construction of the pointed monoid
$\Pi$. This is proved in~\cite[Thm.~7.1]{hm} but see
also~\cite[Prop.~4]{h3}. The cyclic bar-construction is the geometric
realization of the pointed cyclic set with $m$-simplices
$$N^{\cy}(\Pi)[m] = \Pi \wedge \dots \wedge \Pi \hskip 6 mm
\text{($m+1$ factors)}$$
and with the Hochschild-type cyclic structure maps
$$\begin{aligned}
d_i(\pi_0 \wedge \dots \wedge \pi_m) & = \pi_0 \wedge \dots \wedge
\pi_i\pi_{i+1} \wedge \dots \wedge \pi_m, & 0 \leqslant i < m, \cr
{} & = \pi_m\pi_0 \wedge \pi_1 \wedge \dots \wedge \pi_{m-1}, &
i = m, \cr
s_i(\pi_0 \wedge \dots \wedge \pi_m) & = \pi_0 \wedge \dots \wedge
\pi_i \wedge 1 \wedge \pi_{i+1} \wedge \dots \wedge \pi_m, &
0 \leqslant i \leqslant m, \cr
t_m(\pi_0 \wedge \dots \wedge \pi_m) & = \pi_m \wedge \pi_0 \wedge
\pi_1 \wedge \dots \wedge \pi_{m-1}. & \cr
\end{aligned}$$
It is a pointed $\mathbb{T}$-space by the theory of cyclic
sets~\cite[7.1.9]{loday}. It follows that the $\mathbb{T}$-spectrum
$T(A,B,I)$ is canonically $\mathcal{F}$-equivalent to the iterated
mapping fiber of the following diagram of $\mathbb{T}$-spectra.
$$\xymatrix{
{ T(k) \wedge N^{\cy}(\Pi^2) } \ar[r]^{\epsilon} \ar[d]^{(\phi,\phi')} &
{ T(k) \wedge N^{\cy}(\Pi^0) } \ar[d]^{\Delta} \cr
{ (T(k) \wedge N^{\cy}(\Pi^1)) \times (T(k) \wedge N^{\cy}(\Pi^1)) }
\ar[r]^{\epsilon \times \epsilon} &
{ (T(k) \wedge N^{\cy}(\Pi^0)) \times (T(k) \wedge N^{\cy}(\Pi^0)). } \cr
}$$
The cyclic bar-constructions of $\Pi^1$ and $\Pi^2$ have natural
wedge-decompositions which we now explain.

We define $N^{\cy}(\Pi^1,i)[m]$ to be the subset of
$N^{\cy}(\Pi^1)[m]$ that consists of the base-point and of the
simplices $z^{i_0} \wedge \dots \wedge z^{i_m}$ with $i_0 + \dots +
i_m = i$. It is clear that the pointed set $N^{\cy}(\Pi^1)[m]$
decomposes as the wedge-sum of the pointed subsets
$N^{\cy}(\Pi^1,i)[m]$ where $i$ ranges over the non-negative
integers. The cyclic structure maps preserve this decomposition, and
hence, the geometric realization decomposes accordingly as a wedge-sum
of pointed $\mathbb{T}$-spaces
$$N^{\cy}(\Pi^1) = \bigvee N^{\cy}(\Pi^1,i)$$
indexed by the non-negative integers. To state the analogous
wedge-decomposition of $N^{\cy}(\Pi^2)$, we first recall the notion of
a cyclical word.

A word of length $m$ with letters in a set $S$ is a function 
$$\omega \colon \{1,2, \dots, m\} \to S.$$
The action by the cyclic group $C_m$ of order $m$ on the set
$\{1,2,\dots,m\}$ by cyclic permutation induces an action on the set
of words of length $m$ in $S$. A cyclical word of length $m$ with
letters in $S$ is an orbit for the action of $C_m$ on the set of words
of length $m$ in $S$. We write $\bar{\omega}$ for the orbit through
$\omega$. By the period of $\bar{\omega}$, we mean the length of the
orbit $\bar{\omega}$. In particular, the set that consists of the
empty word is a cyclical word of length $0$ and period $1$.

We associate a word $\omega(\pi)$ with letters in $x$ and $y$ to every
non-zero element $\pi \in \Pi^2$. A non-zero element $\pi \in \Pi^2$ is
either of the form $\pi = x^i$ or $\pi = y^i$. In the former case, we
define $\omega(\pi)$ to be the unique word of length $i$ all of whose
letters are $x$, and in the latter case, we define $\omega(\pi)$ to be
the unique word of length $i$ all of whose letters are $y$. More
generally, we associate to every $(m+1)$-tuple $(\pi_0, \dots, \pi_m)$
of non-zero elements of $\Pi$ the word
$$\omega(\pi_0, \dots, \pi_m) = \omega(\pi_0) * \dots *
\omega(\pi_m)$$
defined to be the concatenation of the words $\omega(\pi_0), \dots,
\omega(\pi_m)$. Now, for every cyclical word $\bar{\omega}$ with
letters $x$ and $y$, we define
$$N^{\cy}(\Pi^2,\bar{\omega})[m] \subset N^{\cy}(\Pi^2)[m]$$
to be the subset that consists of the base-point and the elements
$\pi_0 \wedge \dots \wedge \pi_m$, where $(\pi_0, \dots, \pi_m)$
ranges over all $(m+1)$-tuples of non-zero elements of $\Pi^2$ such
that $\omega(\pi_0, \dots, \pi_m) \in \bar{\omega}$. As $m \geqslant 0$
varies, these subsets form a cyclic subset
$$N^{\cy}(\Pi^2,\bar{\omega})[-] \subset N^{\cy}(\Pi^2)[-],$$
and we define $N^{\cy}(\Pi^2,\bar{\omega}) \subset N^{\cy}(\Pi^2)$ to
be the geometric realization. It is clear that the cyclic set
$N^{\cy}(\Pi^2)[-]$ decomposes as the wedge-sum of the cyclic subsets
$N^{\cy}(\Pi^2,\bar{\omega})[-]$, where $\bar{\omega}$ ranges over all
cyclical words with letters $x$ and $y$. Hence, the geometric
realization decomposes as a wedge-sum of pointed $\mathbb{T}$-spaces
$$N^{\cy}(\Pi^2) = \bigvee N^{\cy}(\Pi^2,\bar{\omega})$$
indexed by all cyclical words with letters $x$ and $y$. 

\begin{lemma}\label{birelativethh}There is a canonical
$\mathcal{F}$-equivalence of $\mathbb{T}$-spectra
$$\bigvee T(k) \wedge N^{\cy}(\Pi^2,\bar{\omega}) \xto{\sim}
T(A,B,I),$$
where the wedge-sum on the left-hand side ranges over all cyclical
words of period $s \geqslant 2$ with letters $x$ and $y$.
\end{lemma}

\begin{proof}Let $\bar{\omega}$ be a cyclical word of period $s \geq
2$. Then every word $\omega \in \bar{\omega}$ has both of the letters
$x$ and $y$. Therefore, the compositions of the canonical map
$$T(k) \wedge N^{\cy}(\Pi^2,\bar{\omega}) \to T(A)$$
and the maps $\phi \colon T(A) \to T(B)$, $\phi' \colon T(A) \to
T(B)$, and $\epsilon \colon T(A) \to T(A/I)$ are all equal to the
constant map. Hence, we obtain a canonical map of $\mathbb{T}$-spectra
$$T(k) \wedge N^{\cy}(\Pi^2, \bar{\omega}) \to T(A,B,I)$$
and the wedge-sum of these maps constitute the map of the
statement. The diagram~(\ref{square}) and the 
$\mathcal{F}$-equivalence~(\ref{separate}) shows that this map is an
$\mathcal{F}$-equivalence.
\end{proof}

\begin{lemma}\label{barconstruction}Let $\bar{\omega}$ be a cyclical
word of period $s \geqslant 2$ with letters $x$ and $y$. The homotopy type
of the pointed $\mathbb{T}$-space $N^{\cy}(\Pi^2,\bar{\omega})$ is given
as follows.

(i) If $\bar{\omega}$ has period $s = 2$ and length $m = 2i$, then a
choice of representative word $\omega \in \bar{\omega}$ determines a
$\mathbb{T}$-equivariant homeomorphism
$$S^{\R[C_m] - 1} \wedge_{C_i} \mathbb{T}_+ \xto{\sim}
N^{\cy}(\Pi^2,\bar{\omega}),$$
where $\R[C_m] - 1$ is the reduced regular representation of $C_m$.

(ii) If $\bar{\omega}$ has period $s > 2$, then
$N^{\cy}(\Pi^2,\bar{\omega})$ is $\mathbb{T}$-equivariantly
contractible. 
\end{lemma}

\begin{proof}We refer the reader to~\cite{loday} for the basic theory
of cyclic sets and their geometric realization. Let $\bar{\omega}$ be
a cyclical word of period $s \geq 2$ and length $m = si$ with letters
$x$ and $y$, and let $\omega \in \bar{\omega}$. We let
$(\pi_0,\dots,\pi_{m-1})$ be the unique $m$-tuple of non-zero elements
in $\Pi^2$ such that $\omega(\pi_0,\dots,\pi_{m-1}) = \omega$. Then
the pointed cyclic set $N^{\cy}(\Pi^2,\bar{\omega})[-]$ is generated
by the $(m-1)$-simplex $\pi_0 \wedge \dots \wedge \pi_{m-1}$. Hence,
there is a unique surjective map of pointed cyclic sets
$$f_{\omega} \colon \Lambda^{m-1}[-]_+ \to
N^{\cy}(\Pi^2,\bar{\omega})[-]$$
that maps the canonical generator of the cyclic standard
$(m-1)$-simplex to the generator $\pi_0 \wedge \dots \wedge
\pi_{m-1}$. We recall that the automorphism group of the pointed
cyclic set $\Lambda^{m-1}[-]_+$ is a cyclic group order $m$ generated
by the dual of the cyclic operator $t_m$. Since the cyclic operator $t_m^s$
fixes the generator $\pi_0 \wedge \dots \wedge \pi_{m-1}$, we obtain a
factorization of the map $f_{\omega}$ over the quotient by the
subgroup of the automorphism group or order $i = m/s$,
$$f_{\omega} \colon (\Lambda^{m-1}[-]/C_i)_+ \to
N^{\cy}(\Pi^2,\bar{\omega})[-].$$
We next recall that the geometric realization of the cyclic standard
$(m-1)$-simplex is $\mathbb{T}$-equivariantly homeomorphic to
$\Delta^{m-1} \times \mathbb{T}$, where $\mathbb{T}$ acts by
multiplication in the second factor. Moreover, the homeomorphism may
be chosen in such that the dual of the cyclic operator $t_m$ acts on
$\Delta^{m-1}$ by the affine map that cyclically permutes the vertices
and on $\mathbb{T}$ by rotation through $2\pi/m$;
see~\cite[Sect.~7.2]{hm}. It follows that the map $f_{\omega}$ gives
rise to a continuous $\mathbb{T}$-equivariant surjection
$$f_{\omega} \colon (\Delta^{m-1} \times_{C_i} \mathbb{T})_+ \to
N^{\cy}(\Pi^2,\bar{\omega}).$$
There is a canonical $C_i$-equivariant homeomorphism  
$$\Delta^{s-1} * \dots * \Delta^{s-1} \xto{\sim} \Delta^{m-1},$$
where the group $C_i$ cyclically permutes the $i$ factors in the join
on the left-hand side, and the map $f_{\omega}$ collapses the join of a
number of codimension $1$ faces of $\Delta^{s-1}$ to the base-point.
If the period of $\bar{\omega}$ is $s = 2$, then the map $f_{\omega}$
collapses the whole boundary $\partial\Delta^{m-1} \subset
\Delta^{m-1}$ to the base-point. Hence, in this case, we have a
$\mathbb{T}$-equivariant homeomorphism
$$f_{\omega} \colon (\Delta^{m-1}/\partial \Delta^{m-1}) \wedge_{C_i}
\mathbb{T}_+ \xto{\sim} N^{\cy}(\Pi^2,\bar{\omega}).$$
The simplex $\Delta^{m-1}$ embeds as the convex hull of the
group elements in the regular representation $\R[C_m]$. This
identifies the $C_m$-space $\Delta^{m-1}/\partial \Delta^{m-1}$ with
the one-point compactification of the reduced regular representation
$\R[C_m] - 1$ as stated. This completes the proof of the statement
for $s = 2$. If the period $s > 2$, then there exists a codimension
$1$ face $F \subset \Delta^{s-1}$ that is not collapsed to the
base-point. We have a canonical homeomorphism
$$\operatorname{cone}(F) \xto{\sim} \Delta^{s-1}$$
of the unreduced cone on the face $F$ onto the simplex
$\Delta^{s-1}$. The canonical null-homotopy of the unreduced cone
induces a $C_i$-equivariant null-homotopy
$$\Delta^{s-1} * \dots * \Delta^{s-1} \times [0,1] \to \Delta^{s-1} *
\dots * \Delta^{s-1},$$
and since the face $F$ is not collapsed to the base-point
by the map $f_{\omega}$, this induces a $\mathbb{T}$-equivariant
null-homotopy
$$N^{\cy}(\Pi^2,\bar{\omega}) \wedge [0,1]_+ \to
N^{\cy}(\Pi^2,\bar{\omega}).$$
This completes the proof of the statement for $s > 2$.
\end{proof}

\begin{remark}The statement of Lemma~\ref{barconstruction} may be
viewed as a topological refinement of the calculation in~\cite{bach1}
of the Hochschild homology of the pointed monoid ring $\Z(\Pi^2) =
\Z[x,y]/(xy)$. Indeed, for any pointed monoid $\Pi$, the reduced
singular homology groups $\tilde{H}_*(N^{\cy}(\Pi);\Z)$ and the
Hochschild homology groups $\HH_*(\Z(\Pi))$ are canonically isomorphic.
\end{remark}

\begin{proof}[Proof of Prop.~\ref{thh}]It follows from
Lemmas~\ref{birelativethh} and~\ref{barconstruction} that we have an
$\mathcal{F}$-equivalence of $\mathbb{T}$-spectra
$$\bigvee T(k) \wedge S^{\R[C_{2i}]-1} \wedge_{C_i} \mathbb{T}_+
\xto{\sim} T(A,B,I)$$
where the wedge sum ranges over all positive integers $i$. The
equivalence depends on a choice, for every cyclical word
$\bar{\omega}$ with letters $x$ and $y$ of period $2$, of a
representative word $\omega \in \bar{\omega}$. We choose the
representative $\omega = xy \dots xy$. Now, as a representation of the
subgroup $C_i \subset C_{2i}$, the regular representation $\R[C_{2i}]$
is isomorphic to the complex representation $\lambda_i = \C(1) \oplus
\dots \oplus \C(i)$, where $\C(t)$ denotes the representation of
$\mathbb{T}$ on $\C$ through the $t$-fold power map. Hence, a choice
of such an isomorphism determines a $\mathbb{T}$-equivariant
homeomorphism
$$S^{\lambda_i} \wedge_{C_i} \mathbb{T}_+ \xto{\sim} S^{\R[C_{2i}]}
\wedge_{C_i} \mathbb{T}_+.$$
Moreover, since the $C_i$-action on $\lambda_i$ extends to a
$\mathbb{T}$-action, we further have the canonical
$\mathbb{T}$-equivariant homeomorphism
$$S^{\lambda_i} \wedge (\mathbb{T}/C_i)_+ \xto{\sim}
S^{\lambda_i} \wedge_{C_i} \mathbb{T}_+$$
that takes $(w,zC_i)$ to the class of $(z^{-1}w,z)$. The completes the
proof.
\end{proof}

\section{Topological cyclic homology}

In this section, we prove the formula for the bi-relative topological
cyclic homology groups $\TC_q(A,B,I;p)$ that was stated in
Thm.~\ref{tcthm} of the introduction. We derive this formula from the
corresponding formula for topological Hochschild homology that we
proved in Prop.~\ref{thh} above. The argument is 
very similar to the analogous argument in the case of truncated
polynomial algebras~\cite{hm1,hm2,h3}. We refer the reader
to~\cite[3.7]{h3} for the definition of topological cyclic homology.

We have from Prop.~\ref{thh} an $\mathcal{F}$-equivalence of
$\mathbb{T}$-spectra
$$\bigvee_{i \in \N} T(k) \wedge S^{\lambda_i} \wedge
(\mathbb{T}/C_i)_+[1] \xto{\sim} T(A,B,I),$$
and we wish to evaluate the homotopy groups of the $C_{p^{n-1}}$-fixed
point spectra. To this end, we reindex the wedge-sum on the left-hand side
after the $p$-adic valuation of $i \in \N$. The left-hand side is then
rewritten as
$$\begin{aligned}
{} & \bigvee_{d \in \N} T(k) \wedge S^{\lambda_{p^{n-1}d}} \wedge
(\mathbb{T}/C_{p^{n-1}d})_+[1] \cr
{} & \vee \bigvee_{r=1}^{n-1} \bigvee_{d \in I_p} T(k) \wedge
S^{\lambda_{p^{r-1}d}} \wedge (\mathbb{T}/C_{p^{r-1}d})_+[1], \cr
\end{aligned}$$
where $\mathbb{N}$ and $I_p$ are the sets of positive integers and
positive integers that are not divisible by $p$, respectively. Hence,
the $\mathbb{T}$-spectrum $\rho_{p^{n-1}}^*T(A,B,I)^{C_{p^{n-1}}}$ is
equivalent to the wedge-sum
$$\begin{aligned}
{} & \bigvee_{d \in \N} \rho_{p^{n-1}}^*(T(k) \wedge
S^{\lambda_{p^{n-1}d}} \wedge
(\mathbb{T}/C_{p^{n-1}d})_+)^{C_{p^{n-1}}}[1] \cr
{} & \vee \bigvee_{r=1}^{n-1} \bigvee_{d \in I_p}
\rho_{p^{n-r}}^*(\rho_{p^{r-1}}^*(T(k) \wedge S^{\lambda_{p^{r-1}d}}
  \wedge (\mathbb{T}/C_{p^{r-1}d})_+)^{C_{p^{r-1}}})^{C_{p^{n-r}}}[1]. \cr 
\end{aligned}$$
Now, for every $\mathbb{T}$-spectrum $T$, there is a natural
equivalence of $\mathbb{T}$-spectra
$$\rho_{p^m}^*T^{C_{p^m}} \wedge
\rho_{p^m}^*((\mathbb{T}/C_{p^md})_+)^{C_{p^m}} \xto{\sim}
\rho_{p^m}^*(T \wedge (\mathbb{T}/C_{p^md})_+)^{C_{p^m}}$$
and the $p^m$th root defines a $\mathbb{T}$-equivariant homeomorphism
$$(\mathbb{T}/C_d)_+ \xto{\sim}
\rho_{p^m}^*((\mathbb{T}/C_{p^md})_+)^{C_{p^m}}.$$
Hence, we can rewrite the wedge-sum above as follows.
$$\begin{aligned}
{} & \bigvee_{d \in \N} \rho_{p^{n-1}}^*(T(k) \wedge
S^{\lambda_{p^{n-1}d}} )^{C_{p^{n-1}}} \wedge
(\mathbb{T}/C_d)_+[1]\cr
{} & \vee \bigvee_{r=1}^{n-1} \bigvee_{d \in I_p}
\rho_{p^{n-r}}^*(\rho_{p^{r-1}}^*(T(k) \wedge
S^{\lambda_{p^{r-1}d}})^{C_{p^{r-1}}} \wedge
(\mathbb{T}/C_d)_+)^{C_{p^{n-r}}}[1]. \cr  
\end{aligned}$$
We recall from~\cite[Lemma~3.4.1]{hm3} that if $T$ is a
$\mathbb{T}$-spectrum, if $d \in I_p$, and if $\iota \colon \{C_d\}
\to \mathbb{T}/C_d$ is the canonical inclusion, then the map
$$V^m\iota_* + dV^m\iota_* \colon \pi_q(T) \oplus \pi_{q-1}(T) \to
\pi_q(\rho_{p^m}^*(T \wedge (\mathbb{T}/C_d)_+)^{C_{p^m}})$$
is an isomorphism. It follows that the group $\TR_q^n(A,B,I;p)$ is
canonically isomorphic to the sum
\begin{equation}\label{trsum}
\begin{aligned}
{} & \bigoplus_{d \in \N} \big( \TR_{q+1-\lambda_{p^{n-1}d}}^n(k;p)
\oplus \TR_{q-\lambda_{p^{n-1}d}}^n(k;p) \big) \cr
{} & \oplus \bigoplus_{r = 1}^{n-1} \bigoplus_{d \in I_p}
 \big( \TR_{q+1-\lambda_{p^{r-1}d}}^r(k;p)
\oplus \TR_{q-\lambda_{p^{r-1}d}}^r(k;p) \big). \cr
\end{aligned}
\end{equation}
We consider the groups $\TR_q^n(A,B,I;p)$ for varying $n \geqslant 1$ as a
pro-abelian group whose structure map is the Frobenius map
$$F \colon \TR_q^n(A,B,I;p) \to \TR_q^{n-1}(A,B,I;p).$$
The Frobenius map takes the summand with index $d \in \N$ in the top
line of~(\ref{trsum}) for $n$ to the summand with index $pd \in \N$ in
the top line of~(\ref{trsum}) for $n-1$. It takes the summand with
indices $d \in I_p$ and $1 \leqslant r < n-1$ in the bottom line
of~(\ref{trsum}) for $n$ to the summand with the same indices in the
bottom line of~(\ref{trsum}) for $n-1$. Finally, it takes the summand
with indices $d \in I_p$ and $r = n-1$ in the bottom line
of~(\ref{trsum}) for $n$ to the the summand with index $d \in \N$ in
the top line of~(\ref{trsum}) for $n-1$. It follows that the
sub-pro-abelian group of $\TR_q^n(A,B,I;p)$ given by the top line
of~(\ref{trsum}) is Mittag-Leffler zero, since the sum
in~(\ref{trsum}) is finite. Hence, the projection onto the
quotient pro-abelian group of $\TR_q^n(A,B,I;p)$ given by the bottom
line of~(\ref{trsum}) is an isomorphism of pro-abelian groups. The
value of the Frobenius map on the bottom line of~(\ref{trsum}) follows
immediately from the relations $FV = p$ and $FdV = d$. Indeed, the
Frobenius preserves the direct sum decomposition and restricts to the
maps
$$\begin{aligned}
F = p & \colon \TR_{q+1-\lambda}^r(k;p) \to
\TR_{q+1-\lambda}^r(k;p), \cr
F = \id & \colon \TR_{q-\lambda}^r(k;p) \to
\TR_{q-\lambda}^r(k;p), \cr
\end{aligned}$$
respectively, on the first and second summand of the bottom line
of~(\ref{trsum}).

We now assume that the group $\TR_{q-\lambda}^r(k;p)$ is annihilated
by $p^m$, for some $m$. If $k$ is an $\Fp$-algebra, then this group is
annihilated by $p^r$. For a general ring $k$, we must instead consider
the group $\TR_{q-\lambda}^r(k;p,\Z/p^v)$ which is annihilated by
$p^v$. It follows that the iterated Frobenius $F^m$ induces the zero
map from the first term in the bottom line of~(\ref{trsum}) for $m+n$
to the first term in the bottom line of~(\ref{trsum}) for $n$. Hence,
the canonical projection onto the second term on bottom line
of~(\ref{trsum}),
$$\TR_q^n(A,B,I;p) \to \bigoplus_{r = 1}^{n-1} \bigoplus_{d \in I_p}
\TR_{q-\lambda_{p^{r-1}d}}^r(k;p),$$
is an isomorphism of pro-abelian groups. Here, we recall, the
structure map in the limit system on the left-hand side is the
Frobenius map and in the limit system on the right-hand side is the
canonical projection. The group
$\smash{\TR_{q-\lambda_{p^{r-1}d}}^r(k;p)}$ is zero, if $q <
\dim_{\R}(\lambda_d) = 2d$, and hence, the limit group is the product
\begin{equation}\label{tf}
\TF_q(A,B,I;p) \xto{\sim} \prod_{r \in \N} \prod_{d \in I_p}
\TR_{q-\lambda_{p^{r-1}d}}^r(k;p).
\end{equation}
We can now evaluate the bi-relative topological cyclic homology groups
that are given by the long-exact following long-exact sequence.
$$\dots \to \TC_q(A,B,I;p) \to \TF_q(A,B,I;p) \xto{R - \id}
\TF_q(A,B,I;p) \to \dots$$
Indeed, under the isomorphism~(\ref{tf}), the map $R$ corresponds to
the endomorphism of the product on the right-hand side of~(\ref{tf})
that is induced from the map 
$$R \colon \TR_{q-\lambda_{p^{r-1}d}}^r(k;p) \to
\TR_{q-\lambda_{p^{r-2}d}}^{r-1}(k;p).$$
Hence, the kernel of the map $R-\id$ in~(\ref{tf}) is identified with
the limit
$$\prod_{d \in I_p} \lim_R \TR_{q-\lambda_{p^{r-1}d}}^r(k;p)$$
and the cokernel is identified with the corresponding derived limit.
The following Lemma~\ref{mittagleffler} shows, in particular, that the
limit system satisfies the Mittag-Leffler condition. Hence, the
derived limit vanishes and we obtain an isomorphism
$$\TC_q(A,B,I;p) \xto{\sim} \prod_{d \in I_p} \lim_R
\TR_{q-\lambda_{p^{r-1}d}}^r(k;p)$$
as stated in Thm.~\ref{tcthm}.

\begin{lemma}\label{mittagleffler}The restriction map
$$R \colon \TR_{q-\lambda_{p^{r-1}d}}^r(k;p) \to
\TR_{q-\lambda_{p^{r-2}d}}^{r-1}(k;p)$$
is an isomorphism, for $q < 2p^{r-1}d$.
\end{lemma}

\begin{proof}We recall from~\cite[Thm.~2.2]{hm} that there is a
long-exact sequence 
$$\dots \to
\mathbb{H}_q(C_{p^{r-1}},T(k) \wedge S^{\lambda_{p^{r-1}d}}) \to
\TR_{q-\lambda_{p^{r-1}d}}^r(k;p) \xto{R}
\TR_{q-\lambda_{p^{r-2}d}}^{r-1}(k;p) \to \cdots$$
and that the left-hand groups are given by a spectral sequence
$$E_{s,t}^2 = H_s(C_{p^{r-1}},\TR_{q-\lambda_{p^{r-1}d}}^1(k;p))
\Rightarrow
\mathbb{H}_q(C_{p^{r-1}},T(k) \wedge S^{\lambda_{p^{r-1}d}}).$$
The groups in the $E^2$-term do not depend on the representation
$\lambda_{p^{r-1}d}$ beyond its dimension, and they are zero if $t <
\dim_{\mathbb{R}}(\lambda_{p^{r-1}d}) = 2p^{r-1}d$. It follows that
the map $R$ is an isomorphism, if $\smash{q < 2^{p^{r-1}}d}$ as
stated.
\end{proof}

\section{Regular $\Fp$-algebras}

Let $k$ be a regular $\Fp$-algebra. The structure of the groups
$\TR_{q-\lambda}^n(k;p)$ that occur on the right-hand side of the
statement of Thm.~\ref{tcthm} of the introduction is given
by~\cite[Thm.~2.2.2]{hm2}, but see also ~\cite[Thm.~11]{h3}. If
$\lambda$ is a finite dimensional complex $\mathbb{T}$-representation,
we define
$$\ell_r = \ell_r(\lambda) = \dim_{\C}(\lambda^{C_{p^r}})$$
and $\ell_{-1} = \infty$ such that we have a descending sequence
$$\infty = \ell_{-1} \geqslant \ell_0 \geqslant \ell_1 \geqslant \dots \geqslant \ell_r
\geqslant \ell_{r+1} \geqslant \dots \geqslant \ell_{\infty} =
\dim_{\C}(\lambda^{\mathbb{T}}).$$
Then the following result is \cite[Thm.~2.2.2]{hm2}.

\begin{theorem}\label{ROG}Let $k$ be a regular\, $\Fp$-algebra, and let
$\lambda$ be a finite dimensional complex $\mathbb{T}$-representation. 
There is a canonical isomorphism of abelian groups
$$\bigoplus W_s\Omega_k^{q -2m} \xto{\sim}
\TR_{q-\lambda}^n(k;p)$$
where the sum runs over all integers $m \geqslant \ell_{\infty}$, and where
$s = s(n,m,\lambda)$ is the unique integer such that $\ell_{n-s} \leq
m < \ell_{n-1-s}$. The group $W_s\Omega_k^j$ is understood to be 
zero for non-positive integers $s$. \hfill\space\qed
\end{theorem}

\begin{remark}It appears to be an important problem to determine
the structure of the $RO(\mathbb{T})$-graded equivariant homotopy
groups 
$$\TR_{\alpha}^n(k;p) =
[S^{\alpha} \wedge (\mathbb{T}/C_{p^{n-1}})_+, T(k)]_{\mathbb{T}}$$
for a general virtual $\mathbb{T}$-representation $\alpha$. The
precise definition of $RO(G)$-graded equivariant homotopy groups in
given in~\cite[Appendix]{lewismandell}. One might well hope that the
$RO(\mathbb{T})$-graded equivariant homotopy groups
$\TR_{\alpha}^n(k;p)$ admit an algebraic description similar to that
of the $\Z$-graded equivariant homotopy groups $\TR_q^n(k;p)$ given
in~\cite[Thm.~B]{h}.
\end{remark}

We can now complete the proof of Thm.~\ref{main} of the introduction.

\begin{proof}[Proof of Thm.~\ref{main}]It suffices by Thm.~\ref{tcthm}
to show that for positive integers $d$ prime to $p$, there is a
canonical isomorphism of abelian groups
$$\bigoplus_{m \geqslant 0} W_s\Omega_k^{q-2m} \xto{\sim}
\lim_R \TR_{q-\lambda_{p^{r-1}d}}^r(k;p)$$
where $s = s(m,d)$ is the unique integer that satisfies $p^{s-1}d \leq
m < p^sd$. It follows from~\cite[Thm.~1.2]{hm} that the canonical
projection
$$\lim_R \TR_{q-\lambda_{p^{r-1}d}}^r(k;p) \to
\TR_{q - \lambda_{p^{n-1}d}}^n(k;p)$$
is an isomorphism for $q < \dim_{\R}(\lambda_{p^nd}) =
2p^nd$. But for $\lambda = \lambda_{p^nd}$ we have
$$\ell_{-1} = \infty \geqslant \ell_0 = p^nd \geqslant \ell_1 = p^{n-1}d \geq
\dots \geqslant \ell_n = d \geqslant \ell_{n+1} = \ell_{\infty} = 0.$$
Hence, Thm.~\ref{ROG} gives a canonical isomorphism of abelian groups
$$\bigoplus_{m \geqslant 0} W_s\Omega_k^{q-2m} \xto{\sim}
\TR_{q-\lambda_{p^{n-1}d}}^n(k;p)$$
where $s = s(n,m,\lambda)$ is the minimum of $n$ and the unique
positive integer $t$ that satisfies $p^{t-1}d \leqslant m < p^td$. The
statement follows.
\end{proof}

\begin{remark}We conclude this paper with a conjecture on the
relationship of the $K$-groups of the rings $k[x,y]/(xy)$ and
$k[t]/(t^2)$. The element $f = x - y$ of the ring $A = k[x,y]/(xy)$ is
a non-zero-divisor with quotient ring $A/fA = k[t]/(t^2)$. It follows
that as an $A$-module $A/fA$ has projective dimension $1$, and hence
we have a push-forward  map on the associated $K$-groups
$$i_* \colon K_q(k[t]/(t^2)) \to K_q(k[x,y]/(xy)).$$
The additivity theorem implies that the image of the map $i_*$ is
contained in the subgroup $K_q(A,I)$. In particular, we obtain an
induced push-forward map
$$i_* \colon K_q(k[t]/(t^2),(t)) \to K_q(k[x,y]/(xy),(x,y)).$$
For $k$ a regular $\Fp$-algebra, the relative $K$-groups on the right
and left-hand sides were evaluated in Thm.~\ref{main}
and~\cite[Thm.~A]{hm2}, respectively. On the one hand, there is a
natural long-exact sequence of abelian groups
$$\cdots \to K_q(k[t]/(t^2),(t)) \xto{\partial}
\bigoplus_{m \geqslant 1} \mathbf{W}_m\Omega_k^{q-2m} \xto{V_2}
\bigoplus_{m \geqslant 1} \mathbf{W}_{2m}\Omega_k^{q-2m} \to \cdots$$
and on the other hand, there is a canonical isomorphism of
abelian groups
$$I \colon \bigoplus_{m \geqslant 1} \mathbf{W}_m\Omega_k^{q-2m}
\xto{\sim} K_q(k[x,y]/(xy),(x,y)).$$
We conjecture that the composite $I \circ \partial$ is equal to
the push-forward map $i_*$.
\end{remark}

\providecommand{\bysame}{\leavevmode\hbox to3em{\hrulefill}\thinspace}
\providecommand{\MR}{\relax\ifhmode\unskip\space\fi MR }
\providecommand{\MRhref}[2]{%
  \href{http://www.ams.org/mathscinet-getitem?mr=#1}{#2}
}
\providecommand{\href}[2]{#2}

\end{document}